\documentclass[letterpaper,12pt]{article}
\usepackage{color,amssymb,enumerate,float,amsmath}

\usepackage{ifpdf}
\ifpdf
	\usepackage[pdftex]{graphicx}
	\usepackage[pdftex,unicode,implicit]{hyperref}
\usepackage{caption}
\usepackage{subcaption}
\usepackage{cite}
\usepackage{graphicx}
	\hypersetup{
  	pdftitle     = {}, 
  	pdfkeywords  = {},
  	pdfauthor    = {},
  	pdfcreator   = {pdf\LaTeXe\ with package \flqq hyperref\frqq},
  	pdfproducer  = {pdf\LaTeXe\ with package \flqq hyperref\frqq},
  	pdfpagemode  = UseNone,  
  	pdffitwindow = true,  
  	unicode      = true,
  	plainpages   = true,
  	colorlinks   = true,  
  	citecolor    = black,  
  	urlcolor     = blue, 
  	linkcolor    = black
	}

\else

  \usepackage[dvips]{graphicx}

\fi 

\makeatletter
\@addtoreset{equation}{section}
\makeatother

\begin{document}
	
\thispagestyle{empty}

\begin{center}
{\bf \LARGE  Accurate analytic approximation for a fractional differential equation with a modified Bessel function term}
\vspace*{15mm}

{\large Byron Droguett}$^{1,a}$,
{\large Pablo Martin}$^{2,a}$
{\large Eduardo Rojas}$^{3,b}$
\vspace{3ex}
{\large Jorge Olivares }$^{4,c}$

{\it $^a$ Department of Physics, Universidad de Antofagasta, 1240000 Antofagasta, Chile.}
\\
{\it $^b$Department of Mechanical Engineering, Universidad de Antofagasta, 1240000 Antofagasta, Chile}
\\
{\it $^c$Department of Mathematics, Universidad de Antofagasta, 1240000 Antofagasta, Chile}\\

\vspace{3ex}

$^1${\tt byron.droguett@uantof.cl\,,} \hspace{.5em}
$^2${\tt pablo.martin@uantof.cl\,, } 
\hspace{.5em}
$^3${\tt
eduardo.rojas@uantof.cl\,, }
\hspace{.5em}
$^4${\tt jorge.olivares@uantof.cl\,.}
\vspace*{15mm}

{\bf Abstract}
\begin{quotation}{\small\noindent
A new analytical approximation function is proposed to accurately fit the solution of a fractional differential equation of order one-half, whose nonhomogeneous term is defined by a modified Bessel function of the first kind. The exact analytical solution of this equation is expressed as the product of two modified Bessel functions. The approximation is constructed using an extended multipoint quasi-rational method, which simultaneously incorporates the series expansion and the asymptotic behavior of the Bessel function. A key modification is introduced in the structure of the fitting function, allowing it to reproduce two terms of the asymptotic expansion instead of only one, thereby improving accuracy for large arguments. Numerical analysis shows that for representative parameter values, the maximum relative error between the proposed fitting function and the exact solution of the fractional differential equation is approximately \(0.18\%\), demonstrating the high precision achieved with only six fitting parameters.
}
\end{quotation}
\end{center}

\thispagestyle{empty}

\newpage

\section{Introduction}

Fractional differential equations extend the classical framework by allowing derivatives of noninteger order, offering a flexible and powerful approach to describe complex dynamical systems \cite{Podlubny1999, Diethelm2010}. They are particularly effective for modeling processes with memory and nonlocal effects, making them valuable tools across many areas of science and engineering \cite{Mainardi2010, Metzler2000, Feng2024}. For instance, the Caputo derivative has been widely employed to model anomalous diffusion, accurately representing deviations from classical diffusion in porous media, biological tissues, and turbulent flows \cite{Metzler2000, Magin2006}. In fluid and solid mechanics, fractional derivatives improve the modeling of material behavior, effectively characterizing stress relaxation, creep, and delayed elasticity in polymers, biomaterials, and composite structures \cite{Bagley1983, Genovese2022}. Beyond these fields, fractional calculus has also proven useful in control theory, signal processing, and electrochemical systems, providing a robust framework for capturing real-world complexities beyond classical methods \cite{Monje2010, Caponetto2021}.

In recent years, fractional calculus has gained increasing relevance in cosmology, particularly in models incorporating viscosity and delay effects. In viscous cosmology, dissipative processes are introduced through an effective pressure that modifies the Friedmann and continuity equations within Eckart’s or Israel–Stewart’s formulations \cite{silva2019, tamayo2022, cruz2020, lepe2017}. Such models successfully describe both early-time inflation and late-time acceleration without invoking scalar fields \cite{maartens1996, brevik2017}. Furthermore, incorporating time-delay terms in cosmological equations accounts for the finite response time of the gravitational system, capturing the cumulative and non-instantaneous adaptation of cosmic fluids. These delay effects, often linked to nonlocal interactions and memory phenomena in quantum gravity, yield a more realistic depiction of the Universe’s evolution and can reproduce inflationary or accelerated expansion consistent with observations \cite{Choudhury2012, palpal2022, micolta2025}.

The combination of fractional derivatives and time delays leads to the formulation of Fractional Time-Delayed Differential Equations, which naturally describe systems where both memory and retardation play key roles. These equations have also been successfully applied in engineering contexts such as control systems, viscoelasticity, and population dynamics, as they capture how the present state depends on the system’s history \cite{geremew2024, huseynov2021, moghaddam2013}. Their analytical treatment often involves Laplace transforms and Mittag-Leffler functions, which generalize the exponential function and serve as powerful tools for describing the evolution of complex systems \cite{lukasz2022, amin2021, chen2021}. Moreover, the stability and existence of solutions for nonlinear fractional delay systems remain crucial to ensuring consistent dynamical behavior across physical and engineering applications \cite{lazarevic2011, pakzad2013, dabiri2018}.

Several analytical and numerical approaches have been developed to solve fractional differential equations. Analytical techniques such as Laplace, Mellin, and Fourier transforms, as well as the Green’s function method, have been successfully applied to linear cases \cite{Podlubny1999}. On the numerical side, researchers have proposed various schemes based on finite differences \cite{Li2012, Albadarneh2016}, finite elements \cite{Zhao2017, Nedaiasl2021}, and the finite volume method \cite{Liu2014, Wang2015}. In recent years, mesh-free methods have emerged as efficient alternatives due to their adaptability to large deformations, complex geometries, and heterogeneous materials, as well as their ability to couple with other computational approaches \cite{Nabian2017, Safari2025_3, Gu2005}. Various mesh-free frameworks have also been proposed for fractional differential equations in physics and engineering \cite{Shakeel2020, Dehghan2016, Safari2025_1, Safari2025_2}.

Equations of the form \( y'(x) = I_\nu(x) \) naturally arise in diffusion, conduction, and wave propagation problems involving cylindrical symmetry or inhomogeneous source terms. However, many real-world systems exhibit nonlocal behavior and memory effects that cannot be captured by classical derivatives. To address these limitations, fractional methods based on the Caputo derivative provide a more accurate and flexible modeling framework, extending the traditional approach to anomalous and memory-dependent dynamics. In this work, we focus on the fractional generalization of this equation and develop an analytical approach to approximate the modified Bessel function across the entire spatial domain, ensuring a consistent description valid for both small and large argument limits.

Recent progress in special and fractional calculus includes the generalization of gamma, beta, hypergeometric, and Bessel functions, along with integral representations, generating relations, and fractional differential equation solutions derived from them \cite{enes2014, gurmej2018, jain2020, jain2016}. A modified Caputo-type fractional iterative method has also been proposed, analyzed for convergence, and validated in engineering contexts, demonstrating superior accuracy and efficiency compared to classical techniques \cite{shams2023}. These developments have introduced unified integral formulas and generalized relations that encompass numerous known and novel cases.

The Bessel equation and its modified form play a central role in modeling a wide range of problems in science and engineering. They naturally emerge in contexts such as electrodynamics, where they describe wave propagation in cylindrical geometries \cite{Jackson1998}; plasma physics, where they aid in analyzing stability and wave behavior in magnetized plasmas \cite{Hasegawa1975}; heat transfer problems involving cylindrical or spherical symmetry \cite{Carslaw1959}; and chemical engineering, where they appear in diffusion–reaction systems \cite{Petrova2009}. 

The general solution of the modified Bessel equation of the first kind, \( I_{\nu}(x) \), is typically expressed as a power series. However, for large or intermediate values of \( x \), this representation requires many terms to achieve high accuracy, making computations costly in practical applications. To overcome this difficulty, various fitting functions have been developed to approximate Bessel functions \cite{Baker1996, Peker2021}. Yet, the accuracy of such approximations is often restricted to specific regions of the domain. A significant improvement was introduced through the Multi-Point Quasi-Rational Approximation (MPQA) technique \cite{MartinBook}, which provides a globally valid analytic approximation for all positive \( x \).

In this work, we emphasize the MPQA method as a robust and efficient semi-analytical technique for constructing high-accuracy approximations to the modified Bessel function \( I_{\nu}(x) \). Unlike traditional series or empirical fittings that are either computationally expensive or locally accurate, the MPQA combines information from both the power-series expansion and asymptotic limits to achieve uniform precision across the entire spatial domain. Using only six fitting parameters, it significantly reduces computational cost while preserving analytical simplicity. Its performance demonstrates clear advantages in solving fractional differential equations with nonhomogeneous Bessel terms. Additional studies for fixed values of \( \nu \) confirm the robustness and generality of the method \cite{Martin2018A, Martin2017, Martin2018B}.

In this research, we address a fractional nonhomogeneous differential equation involving an arbitrary Caputo derivative, where the source term is a modified Bessel function of the first kind. When the derivative order is \( 1/2 \), the solution involves the product of two modified Bessel functions. By applying the MPQA technique, we approximate the exact solution and introduce a highly precise analytic approximation for \( I_\nu \) valid for all \( \nu \in [0,1] \), using only six fitting parameters. The new approximation incorporates additional information from the asymptotic expansion of the Bessel function, combining two hyperbolic and rational terms to capture the first two leading asymptotic contributions. This formulation significantly improves accuracy, especially for intermediate \( x \), and represents a substantial advancement over existing methods \cite{Martin2018A, Martin2017, Martin2018B, Martin2020A, Martin2021A, Martin2022A}.

This paper is organized as follows. Section~2 presents the fractional differential equation under study. Section~3 applies the multi-point quasi-rational approximation technique. Section~4 develops the optimization procedure, and Section~5 concludes the work.

\section{Fractional differential equation and analytic solution}

We consider the nonhomogeneous fractional differential equation
\begin{equation}
    {}^{C}D^\alpha_x y(x) = I_{\nu}(x)\,,
    \label{EqCaputo}
\end{equation}
where the nonhomogeneous term \( I_{\nu}(x) \) denotes the modified Bessel function of the first kind with arbitrary order \( \nu \), and \( \alpha \) represents the order of the Caputo fractional derivative. This type of equation is particularly relevant in the study of fractional calculus because the Caputo derivative maintains compatibility with the classical physical interpretation of boundary and initial conditions, making it suitable for modeling realistic dynamical systems.

The solution of Eq.~(\ref{EqCaputo}) can be obtained using the Laplace transform method. The Laplace transform of the Caputo fractional derivative and of the modified Bessel function are, respectively,
\begin{eqnarray}
\label{LD}
\mathcal{L}\{{}^{C}D^\alpha_x y(x)\} &=& s^\alpha \tilde{y}(s) - \sum_{k=0}^{\lceil \alpha \rceil - 1} s^{\alpha-k-1} y^{(k)}(0)\,, \\
\mathcal{L}\{I_\nu(x)\} &=& \frac{(\sqrt{s^2-1}+s)^{-\nu}}{\sqrt{s^2-1}}\,,
\label{LI}
\end{eqnarray}
where \( \tilde{y}(s) \) is the Laplace transform of \( y(x) \) and \( y^{(k)}(0) \) denotes the initial conditions. This formulation highlights the advantage of using the Caputo operator, since it allows one to impose physically meaningful initial data directly on the solution and its derivatives.

Applying the inverse Laplace transform (\ref{LD}) and (\ref{LI}), the general solution can be expressed as
\begin{eqnarray}
    y(x) &=& \frac{\sqrt{\pi}}{2^{\alpha + 2\nu}} \Gamma(1+\nu)\, x^{\alpha} (ix)^{\nu}\,
    {}_2\tilde{F}_3(a_1,a_2; b_1,b_2,b_3; x^2/2)
    + \sum_{k=0}^{[\alpha]-1} \frac{x^k}{\Gamma(1+k)} y^{(k)}(0)\,,
    \nonumber\\
    &&
    \label{solution}
\end{eqnarray}
where ${}_2\tilde{F}_3$ is the regularized generalized hypergeometric function defined by
\begin{eqnarray}
    {}_2\tilde{F}_3(a_1,a_2; b_1,b_2,b_3; x^2/2) &=&
    \frac{{}_2F_3(a_1,a_2; b_1,b_2,b_3; x^2/2)}{\Gamma(b_1)\Gamma(b_2)\Gamma(b_3)}\,,
    \label{RGH}
\end{eqnarray}
and the parameters are related through
\begin{equation}
    a_1=\frac{1+\nu}{2}=a_2-\frac{1}{2}=\frac{b_3}{2}\,, \qquad
    b_1=a_1+\frac{\alpha}{2}=b_2-\frac{1}{2}\,.
    \label{coeficient}
\end{equation}

For the particular case in which the order of the Caputo derivative is \(\alpha = 1/2\), Eqs.~(\ref{solution}-\ref{coeficient}) simplify remarkably, leading to a compact analytic form in terms of the product of two modified Bessel functions:
\begin{eqnarray}
    y(x) &=&
    i^{\nu}\sqrt{\frac{\pi x}{2}}\, I_{\frac{1}{4}(2\nu - 1)}\!\left(\frac{x}{2}\right)
    I_{\frac{1}{4}(2\nu + 1)}\!\left(\frac{x}{2}\right)\,.
    \label{finalsolution}
\end{eqnarray}
By summing the indices of the Bessel functions in the product of the Eq. (\ref{finalsolution}), it can be verified that the overall order of the resulting term is preserved as \( \nu \), ensuring consistency with the nonhomogeneous source term.

Although Eq.~(\ref{finalsolution}) provides an exact analytic form for \( \alpha = 1/2 \), direct numerical evaluation of the modified Bessel functions can become computationally expensive, especially for large or intermediate arguments where their series representations converge slowly. To overcome this difficulty and to construct a more compact yet accurate analytic representation, we employ the Multi-Point Quasi-Rational Approximation technique.
The MPQA method uses information from both the small- and large-argument asymptotic expansions of \( I_\nu(x) \), together with a rational–hyperbolic combination, to reproduce the correct behavior of the function across the entire domain. This approximation enables efficient evaluation of the solution \( y(x) \) while retaining high accuracy, making it particularly suitable for applications in fractional systems where repeated evaluation of special functions is required. In the following section, we develop this approach in detail and discuss the optimization of the fitting parameters that guarantee the precision of the resulting analytic formula.

\section{Approximation procedure}

The modified Bessel function of the first kind, \( I_\nu(x) \), satisfies the modified Bessel differential equation and plays a crucial role in modeling diffusion, conduction, and wave propagation processes. Its analytical representation is well known both in terms of a power series and an asymptotic expansion \cite{Abramowitz}:
\begin{eqnarray}
    I_\nu(x) &=& \left(\frac{x}{2}\right)^\nu
    \sum_{k=0}^\infty \frac{1}{k!\,\Gamma(k+\nu+1)}
    \left(\frac{x}{2}\right)^{2k}\,,
    \label{seriesI}
    \\
    I_\nu(x) &\sim& \frac{e^{x}}{\sqrt{2\pi x}}
    \left(1 - \frac{4\nu^2-1}{8x}
    + \frac{(4\nu^2-1)(4\nu^2-9)}{2!(8x)^2} - \cdots \right)\,.
    \label{eq_asymp_series}
\end{eqnarray}
While Eq.~(\ref{seriesI}) converges rapidly for small arguments, its computational cost grows for intermediate and large values of \(x\), where the asymptotic form (\ref{eq_asymp_series}) becomes more efficient but less accurate near the origin. To obtain a uniform representation valid across the entire domain, we construct an analytic approximation that smoothly interpolates between both regimes.
In the fractional differential equation introduced in the previous section, the modified Bessel function \( I_\nu(x) \) appears as the nonhomogeneous source term. Numerical evaluation of the corresponding fractional solution often requires repeated computation of \( I_\nu(x) \) and its products for different arguments. Therefore, having a compact analytic approximation of \( I_\nu(x) \) with high accuracy and low computational cost becomes essential for the efficient evaluation of the fractional solution \( y(x) \).  
To this end, we extend the multi-point quasi-rational approximation method, originally proposed in \cite{Martin2017}, by incorporating additional asymptotic information that enhances accuracy for intermediate argument values. This approach captures the correct asymptotic behavior for large \(x\) while maintaining agreement with the series expansion at small \(x\).

The proposed approximation is expressed as a rational–hyperbolic combination weighted by polynomial and fractional terms:
\begin{equation}
    \widetilde{I}_{\nu}(x)=
    \left(\frac{x}{2}\right)^\nu
    \frac{(p_0+p_2x^2)\cosh(x)+(p_1+p_3x^2)\frac{\sinh(x)}{x}}
    {\Gamma(\nu+1)\,(1+\lambda^2x^2)^{\nu/2+1/4}\,(1+q x^2)}\,,
    \label{eq_fit_func}
\end{equation}
where \(p_i\), \(q\), and \(\lambda\) are fitting parameters. The coefficients \(p_i\) and \(q\) are obtained analytically, while \(\lambda\) remains a free optimization parameter to minimize the relative error across the entire domain.  
This expression extends the standard MPQA formulation by combining two hyperbolic functions weighted by polynomials, allowing the inclusion of two leading terms of the asymptotic expansion for large \(x\), instead of only one as in \cite{Martin2017}. Consequently, this formulation ensures greater fidelity for intermediate values of \(x\), where most numerical schemes struggle.
To determine the large-\(x\) coefficients, we equate the two leading terms of Eqs.~(\ref{eq_asymp_series}) and (\ref{eq_fit_func}), for fixed \(\nu\):
\begin{eqnarray}
  \frac{e^{x}}{\sqrt{2\pi x}}\left(1-\frac{4\nu^2-1}{8x}\right)
  &\approx&
  \frac{e^{x}}{\sqrt{2\pi x}}\,\frac{1}{f(\lambda,\nu)q}
  \left(p_2+\frac{p_3}{x}\right)\,,
  \label{eq_large_x}
\end{eqnarray}
where
\begin{eqnarray}
f(\lambda,\nu)&=&\frac{2^{1+\nu}\lambda^{\nu+1/2}\Gamma{(1+\nu)}}{\sqrt{2\pi}}\,.
\end{eqnarray}
Matching terms order by order in Eq. (\ref{eq_large_x}) gives
\begin{equation}
p_2 = q\,f(\lambda,\nu)\,, \qquad
p_3 = \frac{1-4\nu^2}{8}\,f(\lambda,\nu)\,q\,.
\end{equation}

To determine the remaining coefficients, we expand around \(x=0\) using the series for the hyperbolic and binomial functions:
\begin{equation}
\sinh(x)=\sum_{n=0}^\infty\frac{x^{2n+1}}{(2n+1)!}\,,\qquad
\cosh(x)=\sum_{n=0}^\infty\frac{x^{2n}}{(2n)!}\,,
\label{hyperbolic}
\end{equation}
\begin{equation}
(1+(\lambda x)^2)^\beta
= \sum_{n=0}^\infty \binom{\beta}{n} (\lambda x)^{2n}, \qquad |x|<1\,,
\label{binomial}
\end{equation}
with \(\beta=\frac{1}{2}(\nu+\tfrac{1}{2})\) and
\begin{eqnarray}
    \binom{\beta}{n}=\frac{1}{n!}\,\beta(\beta-1)\cdots(\beta-(n-1))\,.
    \label{factorial}
\end{eqnarray}
Let us compare the following terms in  Eq. (\ref{seriesI}) 
\begin{eqnarray}
    I_\nu(x)&=&\frac{1}{\Gamma(1+\nu)}\left(\frac{x}{2}\right)^\nu
    \left(
    1+\frac{1}{1+\nu}\left(\frac{x}{2}\right)^2
    +\frac{1}{2(1+\nu)(2+\nu)}\left(\frac{x}{2}\right)^4
\cdots   
\right)
\,,
\nonumber
\\
\label{seriestruncated}
\end{eqnarray}
with the approximation function (\ref{eq_fit_func}) using Eqs. (\ref{hyperbolic}-\ref{factorial}) 
\begin{eqnarray}
&&
\sum_{n=0}^\infty\left[\frac{(p_0+p_2x^2)}{(2n)!}+\frac{(p_1+p_3x^2)}{(2n+1)!}\right]x^{2n}
=
(1+qx^2)
\nonumber
\\
&&
\times\left(
    1+\frac{1}{1+\nu}\left(\frac{x}{2}\right)^2
    +\frac{1}{2(1+\nu)(2+\nu)}\left(\frac{x}{2}\right)^4
\cdots   
\right) \sum_{n=0}^\infty\binom{\beta}{n}\lambda^{2n}x^{2n}
\,.
\label{masterequation}
\nonumber
\\
\end{eqnarray}
Expanding both the exact series (\ref{seriestruncated})  and the approximation (\ref{masterequation}) and matching coefficients up to fourth order yields the following relations:
\begin{eqnarray}
\label{ecc11}
p_0+p_1 &=& 1\,,\\ 
\label{ecc2}
\frac{p_0}{2}+\frac{p_1}{6}+p_2+p_3 &=& 
q+\beta\lambda^2+\frac{1}{4(1+\nu)}\,,\\
\label{ecc3}
\frac{p_0}{24}+\frac{p_1}{120}+\frac{p_2}{2}+\frac{p_3}{6} &=&
\frac{1}{2} (\beta -1) \beta  \lambda ^4+q \left(\beta  \lambda ^2+\frac{1}{4 (v+1)}\right)
\nonumber\\
&&
+\frac{\beta  \lambda ^2}{4 (v+1)}+\frac{1}{32 (v+1) (v+2)}\,,
\\
\frac{p_2}{24}+\frac{p_3}{120}
&=&
q \left(\frac{1}{2} (\beta -1) \beta  \lambda ^4+\frac{\beta  \lambda ^2}{4 (v+1)}+\frac{1}{32 (v+1) (v+2)}\right)
\nonumber\\
&&
+\frac{(\beta -1) \beta  \lambda ^4}{8 (v+1)}+\frac{\beta  \lambda ^2}{32 (v+1) (v+2)}\,.
\label{ecc55}
\end{eqnarray}
Higher-order terms (\(x^6\) and above) are neglected since they contribute minimally within the target precision of the fit.

\begin{figure}[H]
\centering
\includegraphics[width=0.7\linewidth]{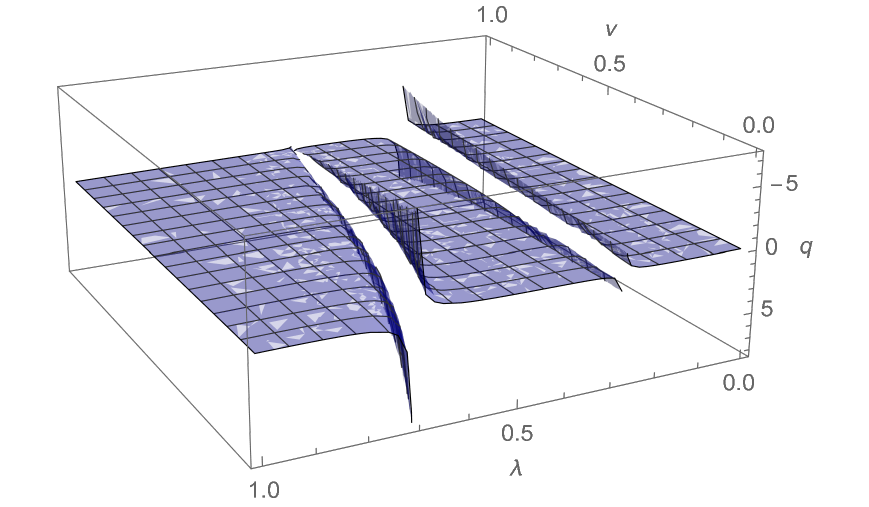}
\caption{Dependence of the parameter \(q\) on \(\lambda\) and \(\nu\). The physically admissible region corresponds to \(q > 0\), ensuring regularity of the approximation function.}
\label{qlv}
\end{figure}

As in Padé-type approximations \cite{Baker1996, Peker2021}, zeros in the denominator (so-called \textit{defects}) can appear, leading to unphysical divergences. To prevent this, both parameters \(q\) and \(\lambda\) are constrained to positive values. Solving the above system Eqs. (\ref{ecc11}-\ref{ecc55}) provides \(q\) as a function of \((\lambda, \nu)\) for all \(\nu \in (0,1)\):
\begin{eqnarray}
    q(\lambda,\nu) &=& 
    \frac{\sqrt{\frac{\pi }{2}}
    \left(-4 \nu^2-240 (\beta -1) \beta  \lambda ^4 (\nu+1) (\nu+2)+24 \beta  \lambda ^2 (\nu+2) (2 \nu-3)+1\right)}
    {4 \left(2^\nu \left(4 \nu^2-49\right) \lambda ^{\nu+\frac{1}{2}} \Gamma (\nu+3)
    +3 \sqrt{2 \pi } (\nu+2) \left(20 \beta  \lambda ^2 (\nu+1)-2 \nu+3\right)\right)}\,.
    \nonumber
    \\
    \label{qnulambda}
\end{eqnarray}
Figure~(\ref{qlv}) shows the dependence of \(q\) on \(\lambda\) and \(\nu\) (see Eq. (\ref{qnulambda})). The region of interest corresponds to the domain where \(q > 0\), ensuring a well-behaved analytic approximation free from defects.
This improved MPQA formulation thus provides a globally accurate analytic expression for \(I_\nu(x)\), combining asymptotic consistency with computational simplicity. In the next section, we discuss the optimization of the free parameter \(\lambda\) and the quantitative analysis of the resulting approximation error.

\section{Optimization}

\begin{figure}[H]
\includegraphics[width=0.7\linewidth]{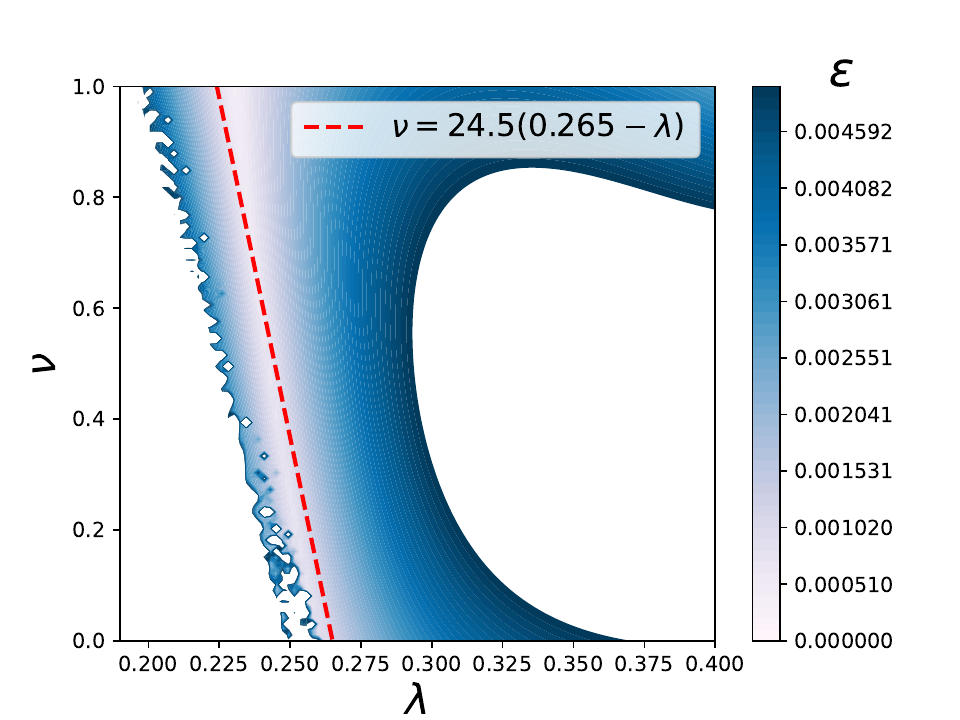}
    \caption{Relative global error (\ref{globalerror}) as a function of the optimization parameter $\lambda$ and $\nu$ for the fitting function}
    \label{fig_error1}
\end{figure}

\begin{figure}[H]
    \centering   
    \includegraphics[width=.48\linewidth]{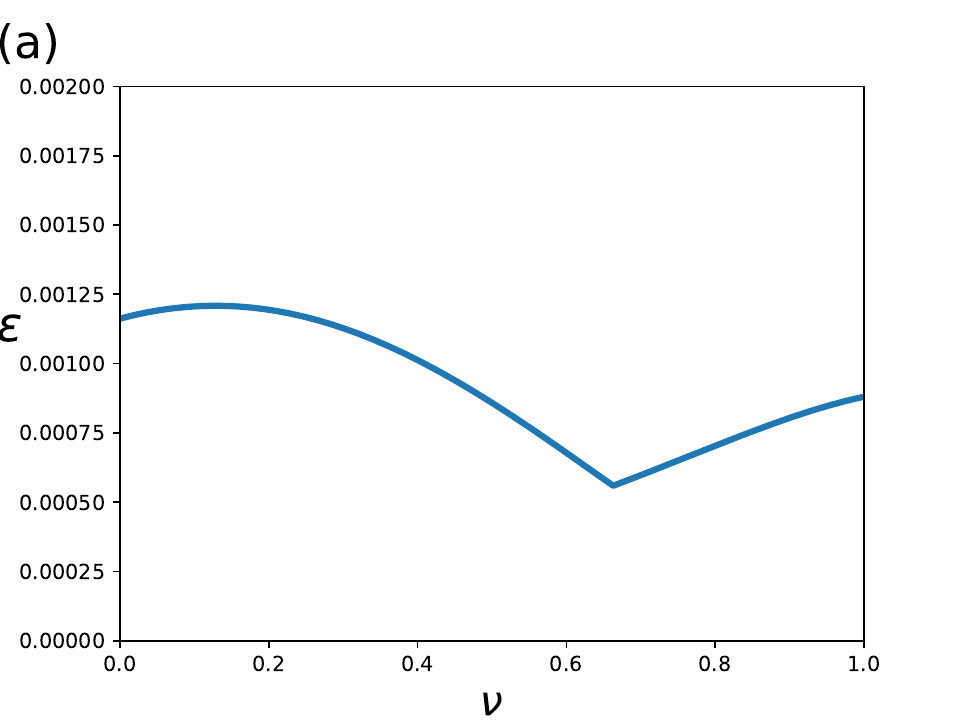} 
    \includegraphics[width=.48\linewidth]{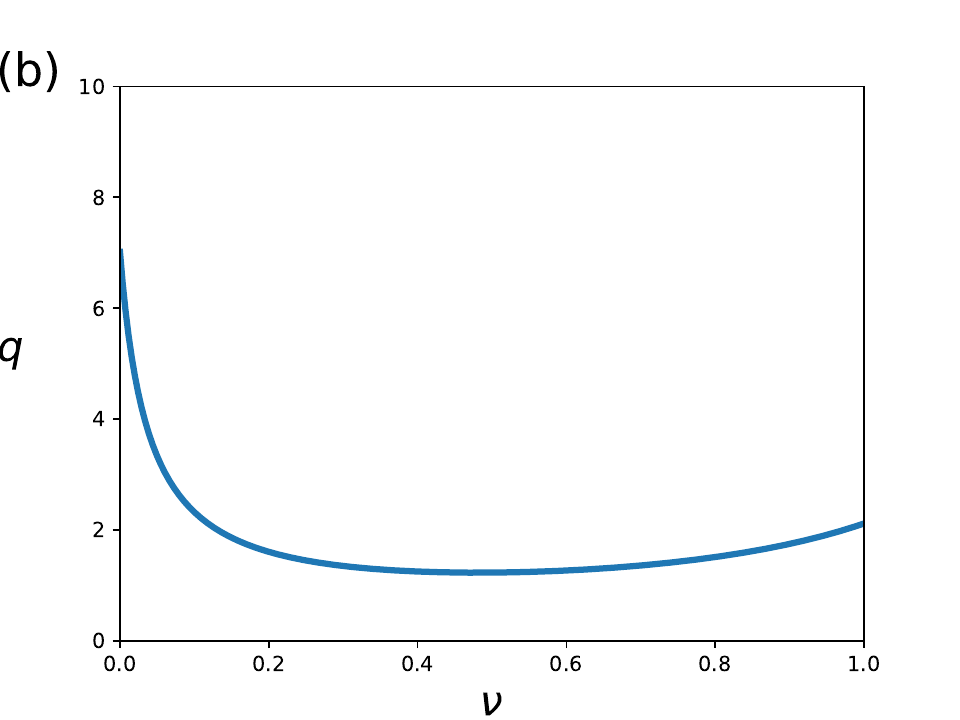}
\caption{ (a) shows the global minimum error as a function of \(\nu\), where for each fixed \(\nu\), the optimal \(\lambda\) parameter is modeled by \(\nu = 24.5(0.265 - \lambda)\). In Fig. (b), we can observe that using this linear approximation, the values of the \(q\) parameter remain positive for all \(\nu\).
}
\label{error3}
\end{figure}

The criterion for determining $\lambda$ is by comparing the maximum error of the fitting function for different values of $\lambda$ and $\nu$. Two types of error definition are established, a punctual error (\ref{punctualerror}) and  a global error (\ref{globalerror}) evaluated on an interval $[a,b]$ of the independent variable $x$. Theses are respectively
\begin{eqnarray}
    \label{punctualerror}
\varepsilon_p(x,\lambda,\nu)&=&
\frac{|I_{\nu}(x)-\widetilde{I}_{\nu}(x,\lambda)|}{I_{\nu}(x)},
\\  \varepsilon(\lambda,\nu)&=&
max\left[{\varepsilon_p(x,\lambda,\nu)}\right]_{x \in[a,b]}.
    \label{globalerror}
\end{eqnarray}
In Fig. (\ref{fig_error1}), the global error (\ref{globalerror}) is shown as a function of the parameters \((\lambda, \nu)\). Several global and relative minima are indicated in the white region of the figure. These values can be modeled by the linear equation \(\nu = 24.5(0.265 - \lambda)\), represented by the red dashed line. Therefore, the optimal values for each parameter of the fitting function can be determined.
In Fig. (\ref{error3}) (a), the global minimum error is shown as a function of the \(\nu\) parameter. We can observe that the error for all values of \(\nu\) is of the order of \(10^{-3}\). In Fig. (\ref{error3}) (b), the \(q\) parameter is shown as a function of \(\nu\) and is always positive. Therefore, there are no issues with poles in the denominator.

\begin{figure}[H]
\includegraphics[width=23pc]{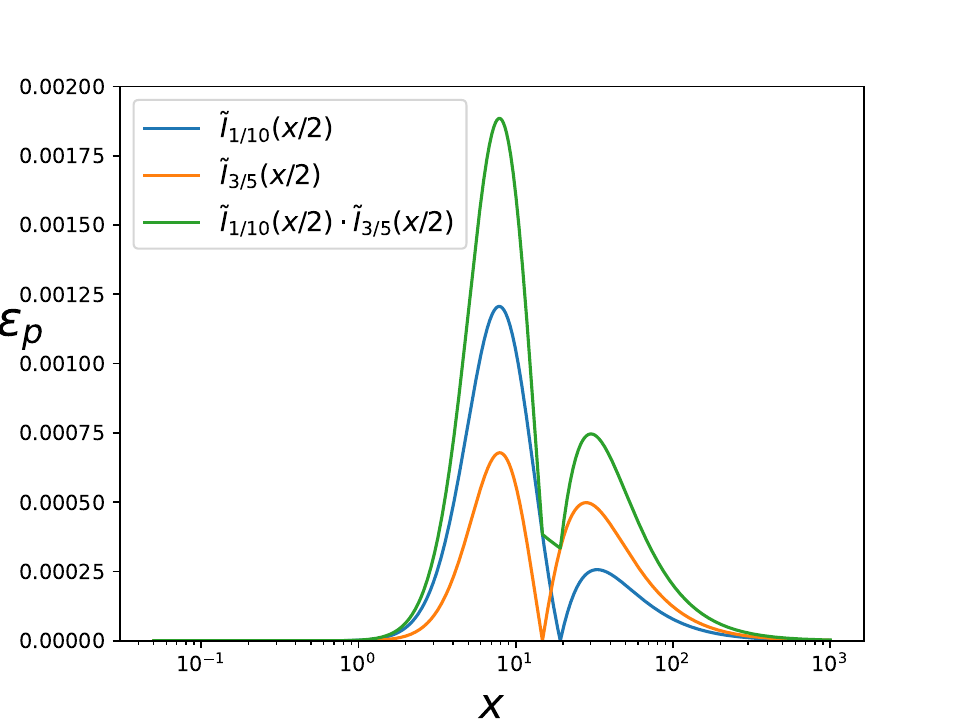}\hspace{2pc} 
\caption{ The relative punctual error (\ref{punctualerror}) as a function of the independent variable \(x\) for the general solution (\ref{finalsolutionexpample}).
}
\label{fig_e_x}
\end{figure}

We present a specific example for the fixed value 
$\nu=7/10$, where Eq. (\ref{finalsolution}) takes the form:
\begin{eqnarray}
    y(x)&=&
    i^{7/10 }\sqrt{\frac{\pi x}{2}}I_{1/10}\left(\frac{x}{2}\right) I_{3/5}\left(\frac{x}{2}\right)\,.
    \label{finalsolutionexpample}
\end{eqnarray}
To apply the approximation method to the solution for this fixed value of \(\nu\), the minimum error \(\varepsilon_{\text{min}}\) is located at \(\lambda_{\text{min}} = 0.236\).
In Fig. (\ref{fig_e_x}), the relative punctual  error (\ref{punctualerror}) is shown for the different analytical fitting functions developed in this study. The maximum errors for \(\widetilde{I}_{1/10}(x)\) and \(\widetilde{I}_{3/5}(x)\) are \(0.1\%\) and \(0.06\%\), respectively. Finally, the maximum relative error for the product of these two modified Bessel functions, \(\widetilde{I}_{1/10}(x)\widetilde{I}_{3/5}(x)\), which corresponds to the error of the solution of the one-half nonhomogeneous fractional differential equation addressed in this work (\ref{finalsolutionexpample}), reaches \(0.18\%\).

Figure (\ref{fig_e_x}) shows that the relative punctual error (\ref{punctualerror}) has a couple of picks for the intermediate values of the argument $x$. One might ask if there are some other picks for large or very large values of $x$. To solve this controversy, a simple and general convergence analysis of the fractional equation solution $y(x) =
    i^{\nu}\sqrt{\frac{\pi x}{2}}\, I_{\nu/2-1/4}\!\left(\frac{x}{2}\right)
    I_{\nu/2+1/4}\!\left(\frac{x}{2}\right)$
    addressed in this study was formulated. Assuming that this exact solution can be expressed by using the asymptotic expansion of each Bessel function for large values of the argument $x/2$, the punctual error is:
\begin{equation}
\begin{aligned}
\varepsilon_p(x,\nu)
&\approx \sum_{k=-\frac{1}{4},\,\frac{1}{4}}
\frac{|I_{\nu/2+k}(x/2)-\widetilde{I}_{\nu/2+k}(x/2,\lambda^*)|}{I_{\nu/2+k}(x/2)} \\[6pt]
&= \sum_{k=-\frac{1}{4},\,\frac{1}{4}}
\frac{\frac{e^{x/2}}{\sqrt{\pi x}}
\left|1-\frac{p_2(\lambda^*,\nu/2+k)}{f(\lambda^*,\nu/2+k)q}
+\left(\frac{1-4(\nu/2+k)^2}{8}
-\frac{p_3(\lambda^*,\nu/2+k)}{f(\lambda^*,\nu/2+k)q}\right)
\frac{2}{x}
+\mathcal{O}\!\left(\frac{1}{x^2}\right)\right|}
{\frac{e^{x/2}}{\sqrt{\pi x}}
\left(1+\frac{1-4(\nu/2+k)^2}{4x}
+\mathcal{O}\!\left(\frac{1}{x^2}\right)\right)} \\[6pt]
&= \sum_{k=-\frac{1}{4},\,\frac{1}{4}}
\left|1-\frac{p_2(\lambda^*,\nu/2+k)}{f(\lambda^*,\nu/2+k)q}\right|
+\mathcal{O}\!\left(\frac{1}{x}\right),
\end{aligned}
\label{eq_error_inf}
\end{equation}
where the parameter $\lambda^*$ is given by the linear approximation of the optimal values of $\lambda$ from Figure (\ref{fig_error1}): 
\begin{equation}
    \lambda^* =\lambda^*(\nu,k) = 0.265-\frac{1}{24.5}\left({\frac{\nu}{2}+k}\right).
\end{equation}
Computing Eq. (\ref{eq_error_inf}):
\begin{eqnarray}
    \sum_{k=-\frac{1}{4},\frac{1}{4}} \left|{1-\frac{p_2(\lambda^*,\nu/2+k)}{f(\lambda^*,\nu/2+k)q}}\right| &\leq&10^{-16}\quad,0\leq\nu/2+k\leq 1
    \,,
\end{eqnarray}
which is of the same order as the truncation error. Therefore, there are no picks for the punctual error when $x$ is large or tends to infinity.

A Python script to create most of the figures presented previously and compute the approximation of the solution for the fractional differential equation presented in this article is available at DOI: \url{https://doi.org/10.6084/m9.figshare.30563453}

\section{Conclusions}

In this work, we obtained the exact analytical solution of a nonhomogeneous fractional differential equation involving a Caputo derivative of arbitrary order \(\alpha\), where the nonhomogeneous term is a modified Bessel function of the first kind. The equation was solved through the Laplace transform, yielding a general solution valid for any order \(\nu\). When the fractional order is fixed to \(\alpha = 1/2\), the solution simplifies elegantly to the product of two modified Bessel functions whose orders sum to that of the nonhomogeneous term, thus preserving the Bessel order structure in the resulting expression.

We then developed a new analytic approximation for the modified Bessel function \(I_\nu(x)\) using the multipoint quasi-rational approximation  method. The proposed function incorporates both the power series and asymptotic expansions of \(I_\nu(x)\), ensuring high accuracy for small, intermediate, and large argument values. A crucial improvement over previous formulations is the inclusion of hyperbolic functions weighted by polynomials, which allows the capture of two leading terms of the asymptotic series instead of one, significantly enhancing precision at intermediate \(x\).
The MPQA formulation depends on six fitting parameters. Four of them are obtained analytically from matching the low-order coefficients of the series expansion, while the remaining two (\(q\) and \(\lambda\)) are determined through an optimization process that minimizes the global relative error. To avoid defects, zeros in the denominator typical of Padé-type approximations the parameters are restricted to positive values, guaranteeing stability of the fitting function over the entire domain. Moreover, the function \(q(\lambda, \nu)\) exhibits a continuous region of positivity for all \(\nu \in [0,1]\), ensuring that the approximation remains well-behaved.

The optimization process revealed an almost linear dependence between the minimal global error and the free parameter \(\lambda\), indicating a smooth and predictable behavior of the fitting accuracy. For representative examples, such as \(\nu = 7/10\), the method yields maximum relative errors of \(0.1\%\) and \(0.06\%\) for \(I_{1/10}\) and \(I_{3/5}\), respectively, and a global error of only \(0.18\%\) for the full solution of the fractional differential equation. Considering that the present method employs only six fitting parameters, these results confirm that the MPQA approach provides an exceptionally accurate, stable, and computationally efficient analytic approximation valid for all \(\nu \in [0,1]\), being more general, accurate, and with fewer parameters than similar approximation techniques like \cite{Salahat2013}. Finally, the strategy can be readily extended to approximate other special functions arising in fractional nonhomogeneous differential equations.

\section*{Acknowledgments}

E. R. thanks Agencia Nacional de Investigación y Desarrollo (ANID-Chile) for financially supporting this research through
Fondecyt Grant 11230970.

\section*{Conflict of Interest}
This work does not have any conflicts of interest.

\end{document}